\documentclass[12pt]{article}
\usepackage{latexsym,amsmath,amssymb}
\newtheorem{theorem}{Theorem}
\newtheorem{corollary}[theorem]{Corollary}

\newtheorem{lemma}[theorem]{Lemma}

\def \R{{\Bbb R}}

\def \remark{\par\noindent{\bf Remark.}\hskip 5pt}
\def \proof{\vskip 10pt \par\noindent{\it Proof.} \hskip 2pt}
\def \proofof#1{\vskip 10pt \par\noindent{\it Proof of #1.} \hskip 2pt}

\def \qed{\hfill $\Box$ \vskip 5pt}
\def \aed{\hfill $\Diamond$ \vskip 5pt}
\def \ve{\varepsilon}
\begin{document}

\title{Goodness of fit test for small diffusions by discrete observations}
\author{Ilia Negri \\
{\small\em
University of Bergamo}\\
{\small\em
Viale Marconi, 5, 24044, Dalmine (BG), Italy}\\
{\small\tt ilia.negri@unibg.it} 
\\ and \\
Yoichi Nishiyama
\footnote{Corresponding author.}
\\
{\small\em
The Institute of Statistical Mathematics} \\
{\small\em
4-6-7 Minami-Azabu, Minato-ku, Tokyo 106-8569, Japan} \\
{\small\tt nisiyama@ism.ac.jp}
}
%\date{November 2007}
\maketitle

\begin{center}
Running title: Goodness of fit test for diffusions
\end{center}

\begin{abstract}
We consider a nonparametric goodness of fit test problem for 
the drift coefficient of one-dimensional small diffusions. 
Our test is based on discrete observation of the processes, 
and the diffusion coefficient is a nuisance function which 
is estimated in our testing procedure. 
We prove that the limit distribution of our test is 
the supremum of the standard Brownian motion, and thus 
our test is asymptotically distribution free. 
We also show that our test is consistent under any 
fixed alternatives. 
\end{abstract}

\begin{center}
\end{center}

\newpage

\section{Introduction}%============================================
Goodness of fit tests play an important role in theoretical
and applied statistics, and the study for them has a long
history.
Such tests are really useful especially if they are
{\em distribution free}, in the sense that their distributions
do not depend on the underlying model.
The origin goes back to the Kolmogorov-Smirnov and
Cr\'amer-von Mises tests in the i.i.d. case,
established early in the 20th century,
and they are
{\em asymptotically distribution free}.
On the other hand, the diffusion process models have been
paid much attention because they are useful in
many applications such as Biology, Medicine, Physics and
Financial Mathematics.
However, the problem of goodness of fit tests for diffusion
processes has still been a new issue in recent years.
Kutoyants \cite{Kut-04} considered this problem in his Section 5.4,
but his tests are not asymptotically distribution free.
Dachian and Kutoyants \cite{Dac-K-08} and
Negri and Nishiyama \cite{Neg-N-07}
proposed some asymptotically distribution free tests.
However, all their results are based on {\em continuous time observation}
of the diffusion processes.
The main contribution of the present paper is that
our test is based on {\em discrete time observation}, which is
more realistic in applications.

Consider a one-dimensional stochastic differential equation (SDE)
\begin{equation}\label{SDE}
X_{t}=x_{0}+\int_{0}^{t}S(X_{s})ds + \ve \int_{0}^{t}\sigma(X_{s})dW_{s},
\quad t \in [0,T],
\end{equation}
where $ S $ and $ \sigma $ are functions which satisfy some properties
described in Section 2, and $ t \leadsto W_{t} $ is a standard Wiener
process defined on a stochastic basis
$ (\Omega, {\cal F}, ({\cal F}_{t})_{t \in [0,T]},P) $.
Here, $ T>0 $ is a fixed time.
We consider a case where a unique strong solution $ X $ to this SDE
exists, and we will consider the asymptotic as $ \ve \downarrow 0 $.
Statistical inference for this model based on continuous observation
was studied by Kutoyants \cite{Kut-94}.
As for discrete observation cases, many researchers have treated
the model in some parametric settings;
see e.g.\ S\o rensen and Uchida \cite{Sor-U-03} and references therein.
In this paper, we are interested in nonparametric
goodness of fit test for the drift coefficient
$ S $, while the diffusion coefficient $ \sigma^2 $ is an unknown
nuisance function which we estimate in our testing procedure.
That is, we consider the problem of testing the hypothesis
$ H_0: S=S_{0} $ versus $ H_1: S\not=S_{0} $ for a given
$ S_{0} $. The meaning of the alternatives ``$ S\not= S_{0} $''
will be precisely stated in Section 4.

We consider the following situation.

\vskip 5pt
\par\noindent
{\bf Sampling Scheme.}
The process $ X=\{ X_{t}; t \in [0,T] \} $
is observed at times
$ 0=t_{0}^{\ve}<t_{1}^{\ve}< \cdots < t_{n(\ve)}^{\ve}=T $,
such that $ h_{\ve}=o(\ve^2) $ as $ \ve \downarrow 0 $,
where $ h_{\ve}=\max_{1 \leq i \leq n(\ve)}|t_{i}^{\ve}-t_{i-1}^{\ve}| $.
\aed

We may assume $ \ve \leq 1 $ and  $ h_\ve \leq 1 $ without loss of generality.
We will propose an asymptotically distribution free test
based on this sampling scheme.

The organization of the article is as follows.
In Section 2, we state some conditions for $ (S,\sigma) $
which are assumed throughout this work.
Section 3 gives the main result under the null hypothesis,
assuming the existence of a consistent estimator for the limit variance.
In Section 4, we prove that our test is consistent under
any fixed alternatives, assuming the existence of a consistent estimator
for the limit variance again.
A consistent estimator for the limit variance is explicitly
constructed in Section 5.
The proofs for lemmas and a theorem in Section 5 will be given in Section 6,
with help from the Appendix.

\section{Preliminaries}%============================================
Let us list some conditions for the pair of functions $ (S,\sigma) $.

\vskip 5pt
\par\noindent
{\bf A1.} There exists
a constant $ C>0 $ such that
\[
|S(x)-S(y)| \leq C |x-y|,
\quad
|\sigma(x)-\sigma(y)| \leq C |x-y|.
\]
\par\noindent
{\bf A2.}
$ \sup_{s \in [0,T]}E|X_{s}|^2 < \infty $.
\aed

\vskip 5pt
Under {\bf A1}, the SDE (\ref{SDE}) has a unique
strong solution $ X $,
and notice also that there exists a constant $ C'>0 $ such that
\[
|S(x)| \leq C'(1+|x|),
\quad
|\sigma(x)| \leq C'(1+|x|) .
\]
To see this, just put $ y=0 $.
The constant $ C' $ depends on the values $ S(0) $ and $ \sigma(0) $,
however the constant $ C $ itself depends on the choice of
the functions $ (S, \sigma) $.
So it is convenient to introduce the notation
\[
K_{S,\sigma}=\max \{ C, C' \}.
\]

Let us fix some more notations.
For given $ S $, let us denote by $ x^{S}= \{ x_{t}^{S}; t \in [0,T] \} $
the solution to the ordinary differential equation
\[
\frac{dx_{t}^{S}}{dt}= S(x_{t}^{S})
\quad \mbox{with the initial value }
x_{0}^{S}=x_{0}.
\]

\vskip 5pt
\par\noindent
{\bf A3.} $ \Sigma_{S,\sigma}:= \sqrt{\int_{0}^{T}\sigma(x_{t}^{S})^{2}dt} >0 $.
\aed

\vskip 5pt
Let us close this section with making some conventions.
We denote by $ C[0,T] $ the space of continuous functions on $ [0,T] $,
and by $ \ell^\infty[0,T] $ the space of bounded functions on $ [0,T] $.
We equip both the spaces with the uniform metric.
We denote by ``$ \to^p $'' and ``$ \to^{d} $'' the convergence
in probability and in distribution as $ \ve \downarrow 0 $, respectively.
The notation ``$ \to $'' always means that we take the limit
as $ \ve \downarrow 0 $.

\section{Asymptotically distribution free test}%=====================

Throughout all this section,
we shall suppose that {\bf A1} - {\bf A3} are satisfied
for some $ (S_{0},\sigma) $.

Our test statistics is based on the random field
$ U^{\ve}=\{ U^{\ve}(u); u \in [0,T] \} $ defined by
\[
U^{\ve}(u)=
\ve^{-1}
\sum_{i=1}^{n(\ve)}1_{[0,u]}(t_{i}^\ve)
[X_{t_i^\ve}-X_{t_{i-1}^\ve}-S_{0}(X_{t_{i-1}^\ve})
|t_{i}^\ve-t_{i-1}^{\ve}|].
\]
We will approximate $ U^{\ve} $ by the following random fields
$ V^{\ve}=\{ V^{\ve}(u); u \in [0,T] \} $ and
$ M^{\ve}=\{ M_{u}^{\ve}; u \in [0,T] \} $,
defined respectively by:
\begin{eqnarray*}
V^{\ve}(u)&=&
\ve^{-1}
\sum_{i=1}^{n(\ve)}1_{[0,u]}(t_{i}^\ve)
\int_{t_{i-1}^{\ve}}^{t_{i}^{\ve}}
[dX_{s}-S_{0}(X_{s})ds];
\\
M_{u}^{\ve}&=&
\ve^{-1}
\int_{0}^{u}
[dX_{s}-S_{0}(X_{s})ds].
\end{eqnarray*}
We present some lemmas which will be proved in Section 6.

\begin{lemma}\label{approximation 1}
$ \sup_{u \in [0,T]}|U^{\ve}(u)-V^{\ve}(u)| \to^p 0 $.
\end{lemma}

\begin{lemma}\label{approximation 2}
$ \sup_{u \in [0,T]}|V^\ve(u)-M_{u}^\ve| \to^p 0 $.
\end{lemma}

\begin{lemma}\label{lemma clt}
$ M^{\ve} \to^d G $ in $ C[0,T] $,
where $ G=\{ G(u); u \in [0,T] \} $ is a Brownian motions
with variance
\[
EG(u)^{2}=\int_{-\infty}^{u}\sigma(x_{t}^{S_0})^{2}dt.
\]
\end{lemma}

Combining these lemmas, we obtain the following result.

\begin{theorem}
$ U^{n} \to^d G $ in $ \ell^{\infty}[0,T] $,
where $ G $ is the process appearing in Lemma \ref{lemma clt}.
\end{theorem}

By the continuous mapping theorem, we have the following.

\begin{corollary}
It holds that
\[
\sup_{u \in [0,T]}|U^{\ve}(u)| \to^d
\sup_{t \in [0,\Sigma_{S_0,\sigma}^{2}]}|B_{t}|
=^d \Sigma_{S_0,\sigma}\sup_{t \in [0,1]}|B_{t}|,
\]
where $ t \leadsto B_{t} $ is a standard Brownian motion,
and the notation ``$ =^d $'' means that the distributions are the same.
\end{corollary}

So we have the main result of the paper.
\begin{theorem}
Under $ H_{0}: S=S_{0} $,
suppose that $ \widehat{\Sigma}^{\ve} $ is a consistent estimator
for $ \Sigma_{S_{0},\sigma} $. Then we have
\[
\frac{\sup_{u \in [0,T]}|U^{\ve}(u)|}{\widehat{\Sigma}^{\ve}}
\to^{d} \sup_{t \in [0,1]}|B_{t}|,
\]
where $ t \leadsto B_{t} $ is a standard Brownian motion.
\end{theorem}

The construction of a consistent estimator $ \widehat{\Sigma}^{\ve} $ for
$ \Sigma_{S,\sigma} $
will be discussed in Section 5.

\section{Consistency of the test}%=========================
Let $ S_{0} $ be that in Section 3.
We denote by $ {\cal S} $ the class of functions $ S $
satisfying {\bf A1} - {\bf A3} and
\begin{equation}\label{consistency assumption}
\int_{0}^{u_S}(S(x_{t}^{S})-S_{0}(x_{t}^{S}))dt
\not= 0
\quad \mbox{for some }u_S \in [0,T].
\end{equation}
The precise description of our problem is testing
the null hypothesis $ H_{0}: S=S_{0} $ versus the alternatives
$ H_{1}: S \in {\cal S} $.

We will prove that our test is consistent.
Fix $ S \in {\cal S} $.
We can write $ U^{\ve}=U_{S}^{\ve}+U_{\Delta}^{\ve} $ where
\[
U_{S}^{\ve}(u)=
\ve^{-1}
\sum_{i=1}^{n(\ve)}1_{[0,u]}(t_{i})[X_{t_i^\ve}-X_{t_{i-1}^\ve}-S(X_{t_{i-1}^\ve})
|t_{i}^\ve-t_{i-1}^{\ve}|]
\]
and
\[
U_{\Delta}^{\ve}(u)=
\ve^{-1}
\sum_{i=1}^{n(\ve)}1_{[0,u]}(t_{i}^\ve)
(S(X_{t_{i-1}^\ve})-S_{0}(X_{t_{i-1}^\ve}))
|t_{i}^\ve-t_{i-1}^{\ve}|.
\]
Now we have
\[
\sup_{u \in [0,T]}
|U^{\ve}(u)|
\geq
\sup_{u \in [0,T]}
|U_{\Delta}^{\ve}(u)|
-
\sup_{u \in [0,T]}
|U_{S}^{\ve}(u)|.
\]
Since $S$ satisfies {\bf A1} - {\bf A3}, by the same argument
as in Section 3, the random field $ U_{S}^{\ve} $
converges to the corresponding Gaussian random field
with $ S_{0} $ replaced by $ S $.
So the second term of the right hand side is $ O_{P}(1) $.
As for the first term of the right hand side, we have the following claim.

\begin{lemma}\label{alternative lemma}
Choose $ u_S \in [0,T] $ as in
(\ref{consistency assumption}).
Then it holds that $ |U_{\Delta}^\ve(u_S)| \not= O_{P}(1) $.
\end{lemma}

We therefore obtain the consistency of the test.
\begin{theorem}
Suppose that $ \widehat{\Sigma}^\ve $ is a consistent estimator for
$ \Sigma_{S,\sigma} $.
Under $ H_{1}: S \in {\cal S} $, it holds that
\[
\frac{\sup_{u \in [0,T]}|U^{\ve}(u)|}{\widehat{\Sigma}^{\ve}}
\not= O_{P}(1).
\]
\end{theorem}

\section{Consistent estimator for $ \Sigma_{S,\sigma} $}

In order to construct an asymptotically distribution free test,
we need a consistent estimator for $ \Sigma_{S,\sigma} $.
The following result gives us an answer.
\begin{theorem}\label{variance}
For any $ (S,\sigma) $ which satisfies {\bf A1} and
$ \sup_{s \in [0,T]}E|X_{s}|^4 < \infty $ (which is stronger than {\bf A2}),
\[
\widehat{\Sigma}^{\ve}=
\sqrt{
\ve^{-2}
\sum_{i=1}^{n(\ve)}
|X_{t_{i}^\ve}-X_{t_{i-1}^\ve}|^2
}
\]
is a consistent estimator for $ \Sigma_{S,\sigma} $.
\end{theorem}

\section{Proofs}%==================================
\proofof{Lemma \ref{approximation 1}}
Without loss of generality, we may assume that
$ \ve \leq 1 $ and $ h_{\ve} \leq 1 $.
It follows from Lemma \ref{Kessler} that
\begin{eqnarray*}
\lefteqn{
E \left( \sup_{u \in [0,T]}|U^{\ve}(u)-V^{\ve}(u)| \right)
}
\\
&\leq&
\ve^{-1}
E
\sum_{i=1}^{n(\ve)}
\int_{t_{i-1}^{\ve}}^{t_{i}^{\ve}}
|S_{0}(X_{t_{i-1}^\ve})
-S_{0}(X_{s})|ds
\\
&\leq&
\ve^{-1}
\sum_{i=1}^{n(\ve)}
\int_{t_{i-1}^{\ve}}^{t_{i}^{\ve}}
K_{S_{0},\sigma}E|X_{t_{i-1}}-X_{s}|ds
\\
&\leq&
\ve^{-1}
T K_{S_{0},\sigma}C_{1}h_{\ve}^{1/2}
\\
&\to&
0.
\end{eqnarray*}
So we have the assertion of the lemma.

\proofof{Lemma \ref{approximation 2}}
Notice that
\begin{eqnarray*}
M_{u}^{\ve}&=&
\ve^{-1}
\sum_{i=1}^{n}
\int_{t_{i-1}^{\ve}}^{t_{i}^{\ve}}
1_{[0,u]}(s)
[dX_{s}-S_{0}(X_{s})ds]
\\
&=&
V^{\ve}(u)
+
\ve^{-1}
\int_{t_{i-1}^{\ve}}^{u}
[dX_{s}-S_{0}(X_{s})ds]
\quad \forall u \in [t_{i-1}^\ve,t_{i}^{\ve})
\\
&=&
V^{\ve}(u)
+
\int_{t_{i-1}^{\ve}}^{u}
\sigma(X_{s})dW_{s}
\quad \forall u \in [t_{i-1}^\ve,t_{i}^{\ve})
\end{eqnarray*}
and that $ M_{T}^{\ve}=V^{\ve}(T) $.
Now we have
\begin{eqnarray*}
E\left|\sup_{u \in [0,T]}|V^{\ve}(u)-M_{u}^{\ve}| \right|^4
&=&
\sum_{i=1}^{n(\ve)}
E\sup_{u \in [t_{i-1}^\ve,t_{i}^\ve)}
|V^{\ve}(u)-M_{u}^{\ve}|^4
\\
&\leq&
\sum_{i=1}^{n(\ve)}
E\sup_{u \in [t_{i-1}^\ve,t_{i}^\ve]}
\left|
\int_{t_{i-1}^{\ve}}^{u}
\sigma(X_{s})dW_{s}\right|^4.
\end{eqnarray*}
It follows from Burkholder-Davis-Gundy's inequality
(see e.g.\ Theorem 26.12 of Kallenberg \cite{Kal-02})
that, for a constant $ c_{k} $ depending only on $ k=4 $,
the right hand side is bounded by
\begin{eqnarray*}
c_{4}
\sum_{i=1}^{n(\ve)}
E
\left|
\int_{t_{i-1}^\ve}^{t_{i}^{\ve}}
\sigma(X_{s})ds \right|^2
&\leq&
c_{4}
\sum_{i=1}^{n(\ve)}
E \left(
\int_{t_{i-1}^\ve}^{t_{i}^{\ve}}
1 ds
\cdot \int_{t_{i-1}^\ve}^{t_{i}^{\ve}}
\sigma(X_{s})^{2}ds \right)
\\
&\leq&
c_{4}
T
\max_{1 \leq i \leq n(\ve)}
\int_{t_{i-1}^\ve}^{t_{i}^{\ve}}
E\sigma(X_{s})^{2}ds
\\
&\leq&
c_{4}
T h_{\ve} \sup_{s \in [0,T]}E\sigma(X_{s})^2
\\
&\to&
0.
\end{eqnarray*}
The proof is finished.
\qed

\proofof{Lemma \ref{lemma clt}}
When $ S=S_{0} $, it holds that
\[
M_{u}^\ve = \int_{0}^u \sigma(X_{s})dW_{s}
\]
We will apply the central limit theorem for continuous martingales.
\begin{eqnarray*}
\langle M^\ve \rangle_{u}
&=&\int_{0}^{u}\sigma(X_{s})^{2}ds
\\
&=&\int_{0}^{u}(\sigma(X_{s})^{2}-\sigma(x_{s}^{S_{0}})^{2})ds
+
\int_{0}^{u}\sigma(x_{s}^{S_{0}})^{2}ds
\\
&=& (I) + (II).
\end{eqnarray*}
Now, using Lemma \ref{Gronwall}, we have
\begin{eqnarray*}
|(I)| & \leq &
\int_{0}^{u}|\sigma(X_{s})^{2}-\sigma(x_{s}^{S_{0}})^{2}|ds
\\
&=&
\int_{0}^{T}|\sigma(X_{s})-\sigma(x_{s}^{S_{0}})|
|\sigma(X_{s})+\sigma(x_{s}^{S_{0}})|ds
\\
&\leq&
K_{S_{0},\sigma}\sup_{t \in [0,T]}|X_{t}-x_{t}^{S_{0}}|
\cdot
\int_{0}^{T}
|\sigma(X_{s})+\sigma(x_{s}^{S_{0}})|ds
\\
&\leq&
K_{S_{0},\sigma}
\exp(K_{S_{0},\sigma}T) \cdot
\ve \sup_{t \in [0,T]}
\left| \int_{0}^{t}\sigma(X_{s})dW_{s} \right|
\cdot
\int_{0}^{T}
|\sigma(X_{s})+\sigma(x_{s}^{S_{0}})|ds
\\
&=& O_{P}(\ve).
\end{eqnarray*}
So we have
$ \langle M^\ve \rangle_{u} \to^p \int_{0}^{u}\sigma(x_{s}^{S_{0}})^{2}ds $,
and the weak convergence of the process $ u \leadsto M_{u}^\ve $ holds.

\proofof{Lemma \ref{alternative lemma}}
We simply denote $ u = u_S $.
We consider the following random variables:
\begin{eqnarray*}
A_{1}^{\ve}&=&
\ve
U_{\Delta}^{\ve}(u)
\\
&=&
\sum_{i=1}^{n(\ve)}1_{[0,u]}(t_{i}^\ve)
\int_{t_{i-1}^\ve}^{t_{i}^\ve}(S(X_{t_{i-1}^\ve})-S_{0}(X_{t_{i-1}^\ve}))ds;
\\
A_{2}^{\ve}&=&
\sum_{i=1}^{n(\ve)}1_{[0,u]}(t_{i}^\ve)
\int_{t_{i-1}^\ve}^{t_{i}^\ve}(S(X_{s})-S_{0}(X_{s}))ds;
\\
A_{3}^\ve &=&
\sum_{i=1}^{n(\ve)}
\int_{t_{i-1}^\ve}^{t_{i}^\ve}
1_{[0,u]}(s)
(S(X_{s})-S_{0}(X_{s}))ds;
\\
A_{4}&=&
\int_{0}^{u}
(S(x_{s}^S)-S_{0}(x_{s}^S))ds.
\end{eqnarray*}

First, it holds that
\begin{eqnarray*}
E|A_{1}^{\ve}-A_{2}^{\ve}|
&\leq&
\sum_{i=1}^{n(\ve)}
E
\int_{t_{i-1}^\ve}^{t_{i}^\ve}
\left\{
|S(X_{t_{i-1}^\ve})-S(X_{s})|
+|S_{0}(X_{t_{i-1}^\ve})-S_{0}(X_{s})|
\right\} ds
\\
&\leq&
(K_{S,\sigma}+K_{S_{0},\sigma})
\sum_{i=1}^{n(\ve)}
E
\int_{t_{i-1}^\ve}^{t_{i}^\ve}
|X_{t_{i-1}^\ve}-X_{s}|ds
\\
&\leq&
(K_{S,\sigma}+K_{S_{0},\sigma})T C_{1}
h_{\ve}^{1/2}
\\
&\to&
0,
\end{eqnarray*}
where $ C_{1} $ is a constant appearing in Lemma \ref{Kessler}.
So we have $ |A_{1}^{\ve}-A_{2}^{\ve}| \to^p 0 $.

Next,
\begin{eqnarray*}
A_{2}^{\ve}-A_{3}^\ve
&=&
\sum_{i=1}^{n(\ve)}
\int_{t_{i-1}^\ve}^{t_{i}^\ve}
(1_{[0,u]}(t_{i}^\ve)-1_{[0,u]}(s))
\{ S(X_{s})-S_{0}(X_{s}) \} ds
\\
&=&
\int_{t_{i-1}^\ve}^{u}
\{ S(X_{s})-S_{0}(X_{s}) \} ds
\quad \forall u \in [t_{i-1}^\ve,t_{i}^\ve).
\end{eqnarray*}
If $ u=u_S=T $, then $ A_{2}^\ve=A_{3}^\ve $.
Since
\begin{eqnarray*}
\lefteqn{
\max_{1 \leq i \leq n(\ve)}
\int_{t_{i-1}^\ve}^{t_{i}^\ve}
\{ |S(X_{s})|+|S_{0}(X_{s})| \} ds
}
\\
&\leq&
h_{\ve}
\sup_{s \in [0,T]} \{ |S(X_{s})|+|S_{0}(X_{s})| \}
=O_{P}(h_{\ve}),
\end{eqnarray*}
we have $ |A_{2}^\ve-A_{3}^\ve| \to^p 0 $.

Finally, notice that
\[
A_{3}^\ve - A_{4}=
\int_{0}^u
(S(X_{s})-S_{0}(X_{s}))ds
-
\int_{0}^u
(S(x_{s}^S)-S_{0}(x_{s}^S))ds.
\]
It follows from Lemma \ref{Gronwall} that
\begin{eqnarray*}
\int_{0}^u
|S(X_{s})-S(x_{s}^S)|ds
&\leq&
\int_{0}^T
K_{S,\sigma}|X_{s}-x_{s}^{S}|ds
\\
&\leq&
T K_{S,\sigma}
\sup_{t \in [0,T]}
|X_{t}-x_{t}^{S}|
\\
&\leq&
T K_{S,\sigma} \cdot
\exp(K_{S,\sigma}T) \cdot
\ve \sup_{t \in [0,T]}
\left| \int_{0}^{t}\sigma(X_{s})dW_{s} \right|
\\
&=&
O_{P}(\ve).
\end{eqnarray*}
By the same way, it holds that
\[
\int_{0}^u
|S_{0}(X_{s})-S_{0}(x_{s}^S)|ds
\leq
T K_{S_{0},\sigma} \cdot
\exp(K_{S,\sigma}T) \cdot
\ve \sup_{t \in [0,T]}
\left| \int_{0}^{t}\sigma(X_{s})dW_{s} \right|
=O_{P}(\ve).
\]
Thus we have $ |A_{3}^\ve-A_{4}| \to 0 $.

Consequently, we obtain $ A_{1}^{\ve} \to^p A_{4} \not=0 $,
which implies that $ |U_{\Delta}^{\ve}(u_S)| \not= O_{P}(1) $.
\qed

\proofof{Theorem \ref{variance}}
By It\^o's formula, we have
\[
|X_{t_i^\ve}|^2-|X_{t_{i-1}^\ve}|^2
=2\int_{t_{i-1}^\ve}^{t_i^\ve}X_{s}dX_s
+ \ve^2 \int_{t_{i-1}^\ve}^{t_i^\ve}\sigma(X_{s})^{2}ds.
\]
Since
\[
|\widehat{\Sigma}^{\ve}|^2
=\ve^{-2}\sum_{i=1}^{n(\ve)}
\left\{
|X_{t_{i}^{\ve}}|^2-|X_{t_{i-1}^{\ve}}|^2
-2X_{t_{i-1}^\ve}(X_{t_{i}^\ve}-X_{t_{i-1}^\ve})
\right\},
\]
it is enough to show that
\[
\ve^{-2}
\sum_{i=1}^{n(\ve)}
\int_{t_{i-1}^\ve}^{t_i^\ve}(X_{s}-X_{t_{i-1}^\ve})dX_s
\to^p 0
\]
and
\[
\int_{0}^{T}
\sigma(X_{s})^{2}ds
\to^p \Sigma_{S,\sigma}^{2}.
\]
The latter is proved by the same argument as that
in the proof of Lemma \ref{lemma clt}.
As for the former, observe that
\begin{eqnarray*}
\lefteqn{
\ve^{-2}
\left|
\sum_{i=1}^{n(\ve)}
\int_{t_{i-1}^\ve}^{t_i^\ve}(X_{s}-X_{t_{i-1}^\ve})dX_s
\right|
}
\\
&\leq&
\ve^{-2}
\sum_{i=1}^{n(\ve)}
\int_{t_{i-1}^\ve}^{t_i^\ve}|X_{s}-X_{t_{i-1}^\ve}||S(X_{s})|ds
+
\ve^{-1}
\left|
\sum_{i=1}^{n(\ve)}
\int_{t_{i-1}^\ve}^{t_i^\ve}(X_{s}-X_{t_{i-1}^\ve})\sigma(X_{s})dW_s
\right|
\end{eqnarray*}
By Lemma \ref{Kessler}, the expectation of the first term on the right hand side is
\begin{eqnarray*}
\lefteqn{
\ve^{-2}
\sum_{i=1}^{n(\ve)}
\int_{t_{i-1}^\ve}^{t_i^\ve}E(|X_{s}-X_{t_{i-1}^\ve}||S(X_{s})|)ds
}
\\
&\leq&
\ve^{-2}
\sum_{i=1}^{n(\ve)}
\int_{t_{i-1}^\ve}^{t_i^\ve}\sqrt{E|X_{s}-X_{t_{i-1}^\ve}|^2}
\sqrt{E|S(X_{s})|^{2}}ds
\\
&\leq&
\ve^{-2}
\sum_{i=1}^{n(\ve)}
\int_{t_{i-1}^\ve}^{t_i^\ve}\sqrt{C_{2} \max \{ h_{\ve}^2, \ve^2 h_{\ve}\} }
\sqrt{E|S(X_{s})|^{2}}ds
\\
&\leq&
\ve^{-2}
T
\sqrt{C_{2}} \max \{ h_{\ve}, \ve h_{\ve}^{1/2} \}
\cdot
\sup_{s \in [0,T]}
\sqrt{E|S(X_{s})|^{2}}
\\
&\to&
0.
\end{eqnarray*}
On the other hand, the expectation of the square of the second term on the right hand side
is
\begin{eqnarray*}
\lefteqn{
\ve^{-2}
\sum_{i=1}^{n(\ve)}
E
\int_{t_{i-1}^\ve}^{t_i^\ve}|X_{s}-X_{t_{i-1}^\ve}|^{2}\sigma(X_{s})^2ds
}
\\
&\leq&
\ve^{-2}
\sum_{i=1}^{n(\ve)}
\int_{t_{i-1}^\ve}^{t_i^\ve}\sqrt{E|X_{s}-X_{t_{i-1}^\ve}|^{4}}
\sqrt{E\sigma(X_{s})^4}ds
\\
&\leq&
\ve^{-2}
T \sqrt{C_4 h_{\ve}^2}
\sup_{s \in [0,T]} \sqrt{E\sigma(X_{s})^4}
\\
&\to&
0.
\end{eqnarray*}
This proves the consistency of our estimator.
\qed

\section*{Appendix}%===============================
In the main part of this article, we use the following
inequality which is well known.
\begin{lemma}\label{Gronwall}
For any solution $ X=\{ X_{t}; t \in [0,T] \} $ to the SDE (\ref{SDE}),
it holds that
\[
\sup_{t \in [0,T]}|X_{t}-x_{t}|
\leq \exp(K_{S,\sigma}T)\cdot \ve \sup_{t \in [0,T]}
\left| \int_{0}^{t}\sigma(X_{s})dW_s
\right|.
\]
\end{lemma}
\proof
The proof is a simple application of Gronwall's inequality:
apply Lemma 4.13 of Liptser and Shiryaev \cite{Lip-S-01} for
$ c_{0}=\ve \sup_{t \in [0,T]} | \int_{0}^{t}\sigma(X_s)dW_s| $,
$ c_{1}= K_{S,\sigma} $,
$ c_{2}=0 $,
$ u(t)=|X_t -x_t| $,
and $ v(t)=1 $.
\qed

The following fact is used in the article
many times, so we state it as a lemma here.
\begin{lemma}
Let $ f:\R \to \R $ be a measurable function such that
$ |f(x)| \leq H(1+|x|) $ for some $ H>0 $.
Let $ X= \{ X_{t}; t \in [0,T] \} $ be any stochastic process.
Let $ k> 0 $ and assume $ \sup_{t \in [0,T]}E|X_{t}|^{k} < \infty $.
Then, it holds that
$ \sup_{t \in [0,T]}E|f(X_{t})|^{k} <\infty $.
\end{lemma}
\proof
Since
\[
|x+y|^{k} \leq ||x|+|y||^k \leq | 2\max \{ |x|,|y| \} |^{k} \leq 2^k \{ |x|^{k} + |y|^{k} \},
\]
the lemma is trivial.
\qed

The following lemma is rather well known, but
we give a full proof for references.
\begin{lemma}\label{Kessler}
Let $ X=\{ X_{t} ; t \in [0,T] \} $ be a solution to the SDE (\ref{SDE})
for $ (S,\sigma) $ which satisfies {\bf A1}.
Let $ k > 0 $ and assume $ \sup_{t \in [0,T]}E|X_{t}|^{k \vee 2} < \infty $.
Then, there exists a constant $ C_{k}>0 $, such that
for any $ 0 \leq t \leq t' \leq T $ and any $ \ve>0 $
\[
E|X_{t'}-X_{t}|^{k}
\leq
C_{k} \max \{ |t'-t|^{k}, \varepsilon^k |t'-t|^{k/2} \}.
\]
In particular, if $ |t'-t|\leq 1 $ and $ \varepsilon \leq 1 $, then
\[
E|X_{t'}-X_{t}|^{k}
\leq
C_{k} |t'-t|^{k/2}.
\]
\end{lemma}
\remark
The constant $ C_{k} $ is {\em not} a universal constant depending only on $ k $.
It actually depends on $ S, \sigma, T $. However, it does not depend on
$ t,t',\varepsilon $.
\proof
First we consider the case $ k \geq 2 $.
Notice that
\[
X_{t'}-X_{t}=
\int_{t}^{t'}
S(X_s)ds
+
\ve
\int_{t}^{t'}
\sigma(X_s)dW_s.
\]
It follows from H\"older's inequality that
\[
\left|
\int_{t}^{t'}
S(X_s)ds
\right|^k
\leq
|t'-t|^{k-1}
\int_{t}^{t'}
|S(X_s)|^kds.
\]
Taking the expectation, it holds that
\[
E
\left|
\int_{t}^{t'}
S(X_s)ds
\right|^k
\leq
|t'-t|^{k}
\sup_{s \in [0,T]}E|S(X_{s})|^{k}.
\]

On the other hand, it follows from Burkholder-Davis-Gundy's
inequality that
there exists a constant $ c_{k}>0 $,
depending only on $ k $, such that
\[
E
\left|
\ve
\int_{t}^{t'}
\sigma(X_s)dW_s
\right|^k
\leq
c_{k}\ve^k
E
\left|
\int_{t}^{t'}
\sigma(X_s)^2ds
\right|^{k/2}.
\]
When $ k > 2 $, it follows from H\"older's inequality that
\begin{eqnarray*}
\ve^k
\left|
\int_{t}^{t'}
\sigma(X_s)^2ds
\right|^{k/2}
&\leq&
\ve^k
|t'-t|^{(k/2)-1}
\int_{t}^{t'}
|\sigma(X_{s})|^{k}ds.
\end{eqnarray*}
Taking the expectation, we have
\[
\ve^k
E
\left|
\int_{t}^{t'}
\sigma(X_s)^2ds
\right|^{k/2}
\leq
\ve^k
|t'-t|^{k/2}
\sup_{s \in [0,T]}E|\sigma(X_{s})|^{k}.
\]
When $ k=2 $, we actually have
\begin{eqnarray*}
\ve^2
E
\int_{t}^{t'}
\sigma(X_s)^2ds
&=&
\ve^2
\int_{t}^{t'}
E \sigma(X_s)^2ds
\\
&\leq&
\ve^2
|t'-t|
\sup_{s \in [0,T]}E \sigma(X_s)^2.
\end{eqnarray*}
Thus the proof for the case $ k \geq 2 $ is finished.

For $ k \in (0,2) $, by Jensen's inequality, it holds that
\[
(E|X_{t'}-X_{t}|^k)^{2/k} \leq E|X_{t'}-X_{t}|^2
\leq
C_{2} \max \{ |t'-t|^{2}, \varepsilon^2 |t'-t| \},
\]
thus we obtain the desired inequality.
\qed

%\vskip 20pt

%\par\noindent
%{\bf Acknowledgements.} 

%\vskip 20pt

\end{document}